\documentclass{amsart}
\usepackage{amssymb}

\theoremstyle{plain}
\newtheorem{theorem}{Theorem}
\newtheorem{corollary}{Corollary}

\newtheorem{lemma}{Lemma}
\newtheorem{proposition}{Proposition}
\newtheorem{conjecture}{Conjecture}

\theoremstyle{definition}
\newtheorem{definition}{Definition}
\newtheorem{example}{Example}

\theoremstyle{remark}

\numberwithin{equation}{section}

\begin{document}

\newcommand{\pd}[2]{ \frac{\partial #1}{\partial #2}}
\newcommand{\AffLag}[2]{ AL_{ #1 } \left ( \mathbb R^{2 #2 } \right ) }
\newcommand{\Lag}[2]{ L_{ #1 } \left ( \mathbb R^{2 #2 } \right ) }
\newcommand{\ASp}[1]{ ASp \left ( \mathbb R^{2 #1} \right ) }

\newlength{\tlen}
\newcommand{\fixUp}[1]{ \( \overset{ }{\underset{ }{#1}} \) }

\title[Affine Symplectic Geometry]
      {Lagrangian submanifolds in affine symplectic\\  
       geometry}
\author{Benjamin McKay}
\address{University College Cork \\
         Cork, Ireland} 
\email{b.mckay@ucc.ie}
\date{\today}
\begin{abstract} We uncover the lowest order 
differential invariants of
Lagrangian submanifolds under affine symplectic maps, and
find out what happens when they are constant.
\end{abstract}
\thanks{The author is grateful for
many fruitful discussions with Robert Bryant.
This work was supported in full or in part by a grant from 
the University of South Florida St. Petersburg New Investigator
Research Grant Fund. This support does not necessarily imply
endorsement by the University of research conclusions.
The author is grateful for the hospitality of 
the Institute for Advanced Study; this paper was written 
there in 1996.}

\maketitle

\section{Introduction} 

In an affine space with translation invariant symplectic
form, I can try to get two Lagrangian submanifolds
into high order of contact, using maps that preserve
the affine and symplectic structure. The lowest
order local invariants
that prevent this are described below, equated with
the space of real cubic hypersurfaces in projective
space, modulo projective automorphism. Setting these
invariants to constants, we will see that there is
still remarkable flexibility in the Lagrangian
submanifold.

\section{Generating Functions}

\begin{definition} An \emph{affine symplectic space}
is a real finite dimensional affine space with translation invariant
symplectic form.
\end{definition}
\begin{example}
$\mathbb{R}^{2n}$ with linear coordinates $q^i,p_i$ and
symplectic form $dp_i \wedge dq^i$ is
an affine symplectic space.
\end{example}

I will work entirely in the smooth category, except
as stated explicitly when I use the Cartan-K\"ahler
Theorem. I record a result from symplectic geometry:
\begin{theorem} Suppose 
that $L$ is a Lagrangian
submanifold of an affine symplectic space, and
we choose a point $x \in L$. Then after an affine 
symplectomorphism,
we can identify our affine symplectic space with
$V \oplus V^*$ (where $V$ is a vector space),
with symplectic form
$$(v,\xi),(v',\xi') \mapsto \left < \xi' , v \right >
                          - \left < \xi , v' \right >$$
so that our point becomes $(0,0)$ and the tangent
space $T_xL$ becomes $V \oplus 0$. In such an
identification, there is a function
$S : V \to \mathbb R$ whose germ at the origin is uniquely
determined, so that near $(0,0)$ the Lagrangian
submanifold is identified with
$$\{ (q,p) : p = S'(q) \}$$
and so that $S(0)=0, S'(0)=0, S''(0)=0$. Every such
function $S$ gives rise to a Lagrangian submanifold
in this way.
\end{theorem}
\begin{proof}
Fixing a point of the affine symplectic
space turns it into a symplectic vector space;
for the rest see \cite{AM}, pp. 161, 402.
\end{proof}

If a Lagrangian submanifold is given
as 
$$\{p = S'(q)\} \subset V \otimes V^*$$
then the function $S$ is called a \emph{generating function}
of the Lagrangian submanifold.
The generating function in the theorem above is
determined once the indicated identification is made;
hence it is determined up to choices of such
identification: the group of linear symplectic
transformations of $V \oplus V^*$ which preserve
the subspace $V \oplus 0$. Mere linear algebra
shows that this group is
$$P = \left \{
\begin{pmatrix}
A & 0 \\
0 & \left (A^t\right )^{-1}
\end{pmatrix}
\begin{pmatrix}
I & B \\
0 & I
\end{pmatrix}
\right \}
$$
where $A \in GL(V)$ and $B : V^* \to V$ is symmetric.
Another point of view: we can take $V$ to be the tangent
space $T_xL$ and then the choice of $S$ is determined
by choice of a Lagrangian subspace complementary to
$T_xL$.
 
\begin{lemma} The subgroup
$$\begin{pmatrix}
I & B \\
0 & I
\end{pmatrix}
$$
of $P$ preserves the lowest order term of the Taylor
expansion of the generating function $S$.
\end{lemma}
\begin{proof} If we change coordinates as
$$\begin{pmatrix} Q \\ P \end{pmatrix}
=
\begin{pmatrix} I & B \\
               0 & I 
\end{pmatrix}
\begin{pmatrix} q \\ p \end{pmatrix}
$$
then in the new coordinates, there will be
some new generating function, $T(Q)$. At all
points $(q,p)$ on our Lagrangian submanifold,
we have $p=S'(q)$ and $P=T'(Q)$, and
$$\begin{array}{l}
Q = q + Bp \\
P = p
\end{array}
$$
so $T'(Q) = S'(q) = S'(Q-BP) = S'(Q-BT'(Q))$ or
$$T'(Q) = S'(Q-BT'(Q)).$$
Since $T$ and $S$ both vanish to second order,
the above equation forces them to have
the same lowest order terms in their Taylor
expansions.
\end{proof}
The abelian group
$$\begin{pmatrix}
I & B \\
0 & I
\end{pmatrix}
$$
thus acts on Taylor expansions as
a subgroup of the nilpotent group
$$\begin{pmatrix}
I & 0 & 0 & \dots  \\
* & I & 0 & \dots  \\
* & * & I & \dots  \\
\vdots & \vdots & \vdots & \ddots
\end{pmatrix}
$$
where each row represents an order of Taylor
coefficients.

\begin{lemma} The group
$$\begin{pmatrix}
A & 0 \\
0 & \left (A^t\right )^{-1}
\end{pmatrix}$$
acts on generating functions by the contragredient
representation of $GL(V)$:
$$S(q) \mapsto S(A^{-1}q).$$ 
\end{lemma}
\begin{proof}
$$\frac{\partial}{\partial q}S(A^{-1}q) 
= \left ( A^{-1} \right)^t 
  S'(A^{-1}q)
$$
\end{proof}

\begin{theorem} The space of 2nd order contact classes
of Lagrangian submanifolds in an affine symplectic space
of dimension $2n$ is canonically identified with the
real cubic hypersurfaces in projective space $\mathbb CP^{n-1}$
modulo real projective automorphisms.
\end{theorem}
\begin{proof} We map each pointed Lagrangian submanifold
to the third order terms in the Taylor expansion of its
generating function. From the last two lemmas, this is
well defined.
\end{proof}

\begin{example}{The Clifford torus:}
$$T^n = \{ (q,p) \in \mathbb R^{2n} :
           q_i^2 + p_i^2 = r_i^2 \}$$
where the $r_i$ are some constants. The reader
may check that the generating function
given near $q=0, p=r$ is naturally
identified, by integrating $p=S'(q)$
with
$$S(q) = \sum \frac{r_i^2}2 
           \left (
             \arcsin \frac{q_i}{r_i} + \frac{q_i}{r_i}
                         \sqrt{ 1-\frac{q_i^2}{r_i^2} }
           \right )
$$
and that the Taylor expansion is
$$S(q) = \sum r_i q_i - \sum \frac{r_i}6 q_i^3 + \dots$$
ignoring 5th and higher order terms. But the subtraction of the 
linear term simply translates the Lagrangian submanifold,
so we may ignore it. Hence, up to affine symplectic
transformation, the Clifford torus has generating
function
$$S(q) = \frac{1}{6} \sum q_i^3 + \dots.$$
\end{example}

\begin{example} Given any cubic form $C \in S^3V^*$ we can
take the graph of its differential
$$L = \{ (q,C'(q)) | q \in V \}.$$
It is easy to see that not only does this have $C$
as its cubic form, but that the 2nd order contact
class of $L$ at \emph{every point} is given by $C$.
Therefore for every cubic form, there is at least
one Lagrangian submanifold which has this cubic form
representing its second order contact class at all points.
\end{example}

\section{Singular Behaviour}

Recall that a homogeneous polynomial on a vector space 
is called \emph{singular}
if there exists a (possibly complex) line
of critical points of the polynomial.
Generically among the homogeneous polynomials
of fixed degree, this does not happen (where ``generic'' 
here indicates
the Zariski topology). If a critical line exists for
a homogeneous cubic polynomial then a real critical
line exists by Bezout's theorem.
Generic Lagrangian submanifolds at general points will have
generating functions with nonsingular cubic terms (where
``generic'' here indicates the smooth topology).
Such a critical line for the cubic term of the
generating function, if real, corresponds to a line with
second order contact with the Lagrangian submanifold.
So generally a Lagrangian submanifold has no such lines, i.e.
it turns away from its own tangent space, like a curve
in Euclidean space. By contrast, a surface
in Euclidean space
at a point of negative Gauss curvature, under 
``generic circumstances'', cuts
its own tangent space on a curve with a double point, hence has
two such lines of 2nd order contact: the two tangent lines
to this curve at the double point.

\begin{example} (1) Consider the function
$$S(q) = \frac{1}{3} \left ( R^2 - |q|^2 \right ) ^{3/2}
\left ( |q|^2 - r^2 \right ) ^{3/2}$$
where $r \le |q| \le R, \ q \in \mathbb R^n$. The
reader may check that
$S^1 \times S^{n-1}$ is immersed into
$\mathbb R^{2n}$ as the Lagrangian submanifold
$$\{(q,\pm S'(q)): r \le |q| \le R \}$$
which is invariant under the linear symplectic maps induced
by rotation of the $q$ variable. 

(2) There is another
way to get $S^1 \times S^{n-1}$ to immerse as 
a Lagrangian submanifold: take any immersion of
$S^1$ into the plane, and any Lagrangian immersion
of $S^{n-1}$ into $\mathbb{R}^{2(n-1)}$, for instance as
$$\{(q,\pm f'(q)) : |q| \le 1\}$$
where
$$f(q) = \frac{1}{3} \left ( 1 - |q|^2 \right ).$$
In either case, by computing the Taylor expansions
of the functions $f$ and $S$ the reader may calculate
that these Lagrangian submanifolds have everywhere
singular cubic forms, but with somewhat different
singularities. (In case (1) there is always a
linear factor to the cubic term of any generating
function about any point, for example).
Moreover these Lagrangian immersions
are distinct under symplectomorphism, if
in case (2) we take the map of the circle to the
plane to be a Jordan curve, for instance, since
in case (1) $\int p \, dq$ over the $S^1$ cycle 
is zero, which in case (2), $\int p \, dq$ is
the area contained inside the Jordan curve.
\end{example}

\section{Cubic Hypersurfaces in Projective Space}
We reach an insurmountable obstacle: cubic
hypersurfaces in projective space can not
be classified up to projective automorphism,
except in low dimensions. Homogeneous polynomials 
of degree d in n variables form a vector space of
$\binom{n+d-1}{d}$ 
dimensions\cite{Hilb}. In the
case of cubic polynomials, there are:
\begin{equation*}
\begin{array}{ll}
  \frac{n(n-1)(n-2)}6 & \mbox{ways to choose distinct indecies } i,j,k \\
  n(n-1) & \mbox{ways to choose } i=j \ne k \\
  n & \mbox{ways to choose } i=j=k
\end{array}  
\end{equation*}
so $\binom{n+2}{3}=\frac{n^3+3n^2+2n}{6}$ dimensions. Since
we have a group of $n^2$ dimensions, acting in a faithful
representation, any sort of reasonable
quotient space must have
$$\binom{n+2}{3} - n^2 = \binom{n}{3}$$
dimensions. For example:
\begin{equation*}
  \begin{array}{ll}
      n & \binom{n}{3} \\ \hline
      3 & 1 \\
      4 & 4 \\
      5 & 10 \\
      6 & 20 \\
      7 & 35
  \end{array}
\end{equation*}
To convince oneself that the group action is
faithful, consider the equivariant map
$$V^* \to S^3V^*, \ \lambda \mapsto \lambda^3.$$

There are three types of cubic forms which we
wish to single out. First, there are the singular ones.
\begin{definition}
Take $V$ a finite dimensional real vector space.
We say that a form $F \in S^3V^*$ is \emph{singular} if 
its zero locus in $\mathbb{CP}V$ has
a singular point, or equivalently if $F$ has a critical
point away from the origin in $V \otimes \mathbb C$.
\end{definition}
\begin{definition}
For forms of any fixed
degree and number of variables, there is an
irreducible homogeneous
polynomial in the coefficients of these forms, called
the \emph{discriminant}, which vanishes precisely on the singular
forms. It is unique up to scaling.
\end{definition}
For cubic forms in $n$ variables, the discriminant
has degree $n2^{n-1}$ (\cite{GKZ}, p.79). Thus the 
discriminant changes
only by a positive factor under linear coordinate
changes.

Secondly, there are the semistable cubic forms (see
\cite{Kir}, pg. 102,  \cite{MFK}, pg. 194).
\begin{definition} A form $F \in S^3V^*$ 
is called \emph{semistable} if there is a
polynomial on $S^3V^*$ (i.e. a polynomial in the 
coefficients of cubic forms), invariant up to scaling
under the action of $GL(V)$ on $S^3V^*$,
which does not vanish on $F$. Semistable forms
constitute a Zariski open set. Forms which are
not semistable are called \emph{null}.
\end{definition}

Thirdly, there are the stable cubic forms. 
\begin{definition} A form $F \in S^3V^*$
is \emph{stable} if
$(1)$ only finitely many linear transformations of $V$
leave $F$ invariant, and $(2)$
there is some polynomial $p$ on
$S^3V^*$ invariant up to scaling, not vanishing
on $F$, so that the orbits of $GL(V)$ on which
$p \ne 0$ are algebraically closed. Stable
forms also constitute a Zariski open set.
\end{definition}

We quote some theorems of geometric invariant theory:

\begin{theorem} Nonsingular $\subset$ stable $\subset$ semistable,
null $=$ singular
\end{theorem}
\begin{proof} See \cite{MFK}, pg. 79.
\end{proof}

\begin{theorem} The standard
representation of the general linear group $GL(V)$ on 
the space of cubic forms $S^3V^*$ has a canonical
choice of moduli space, $\mathfrak M_{3,n}, (n=dim \ V)$, 
given as
the projective variety associated to the algebra
of covariants. There is a well defined regular map
of the semistable elements
$$\left ( S^3V^* \right )^{ss} \to \mathfrak M_{3,n}$$
and the fiber of this map through any stable element
is the orbit of that element.
\end{theorem}

I have a conjecture to make concerning homogeneous
polynomials in general:
\begin{definition} Say that a homogeneous polynomial
$p : \mathbb R^n \to \mathbb R$ of degree $d$ is in \emph{nice} form
if it has the form
$$p(x) = \sum_i (-1)^{\epsilon_i}x_i^d + q(x)$$
where $q(x)$ has no terms of the form $x_i^d$ or $x_i^{d-1}x_j$.
\end{definition}
Alternatively, $p$ is in nice form if it has 1 at each
vertex of its Newton polyhedron, and 0 at each node of
its Newton polyhedron which is next to a vertex.
\begin{definition} A homogeneous polynomial is \emph{nice} if
it can be brought to nice form by a linear coordinate
change.
\end{definition}

\begin{conjecture} Nonsingular implies nice. For degree
greater than 2, there
is precisely one choice of coordinates in which
a nice polynomial reaches nice form, up to
permutation of coordinates, and  (if the degree is even) switching
signs of coordinates.
\end{conjecture}

The result is obvious for degree $d=2$, with any $n$,
and is known for $d=3, n \le 3$. 

\begin{proposition} There is an open set (in the Euclidean
topology) in $S^d \mathbb R^n$ of polynomials which are
nice.
\end{proposition}
\begin{proof} Polynomials in nice form 
are in codimension $n$, and they stay in nice
form only only linear changes of variable
given by matrices whose diagonal terms vanish. 
The result follows by dimension count.
\end{proof}

\begin{corollary} Nice polynomials form
a Zariski constructible set, containing a 
Euclidean open set.
\end{corollary}

Over the complex numbers, this proves that generic
nonsingular homogeneous polynomials are nice.

\begin{corollary} Nice polynomials of odd degree
form a Zariski open set.
\end{corollary}
\begin{proof} Use induction on the number of variables,
by setting each of them to $0$,
and use the fact that reduction to nice form is generically possible 
over complex numbers, to show that the complex linear
transformation reducing a real odd degree polynomial
to nice form must actually be a real linear transformation.
\end{proof}

The hyperbolic geometry of the moduli space
of cubic surfaces (see \cite{ACT}) might be useful in studying
Lagrangian submanifolds.

\section{Three Points on a Projective Line}

\begin{table}
  \begin{tabular}{|l|l|l|} \hline
    \textbf{Normal form}
        & \textbf{Isotropy group} 
        & \textbf{Type} \\ \hline
    $0$                         & \fixUp{GL(2,\mathbb R)} & $0$ \\ \hline
    \( \frac{1}{6}x^3 \)            & \fixUp{  \left \{ 
                                        \begin{pmatrix}
                                        1 & 0 \\
                                        c & d
                                        \end{pmatrix}
                                        : d \ne 0
                                \right \}  }
        & \( \text{linear}^3 \) \\ \hline
    \( \frac{1}{2}x^2y \)       & \fixUp{ \left \{
                                        \begin{pmatrix}
                                        a & 0 \\
                                        0 & \frac{1}{a^2}
                                        \end{pmatrix}
                                        : a \ne 0
                                \right \} }
        &  \( \text{linear}_1 \cdot ( \text{linear}_2 )^2 \) \\ \hline
    \( \frac{1}{6} x^3
+
    \frac{1}{6} y^3 \)          & \fixUp{ \left \{
                                        \begin{pmatrix}
                                        1 & 0 \\ 0 & 1
                                        \end{pmatrix}
                                        ,
                                        \begin{pmatrix}
                                        0 & 1 \\ 1 & 0
                                        \end{pmatrix}
                                  \right \} }
     & (linear)(irred \ quadratic) \\ \hline
     \( \frac{1}{2} x^2 y 
     -
     \frac{1}{2} x y^2 \)         &  \fixUp{\Sigma_3} 
&  \( \text{linear}_1 \cdot \text{linear}_2 \cdot \text{linear}_3 \)
\\ \hline
  \end{tabular}
\caption{Cubic forms in 2 variables}\label{tbl:twovarcubics}
\end{table}
A cubic form in two variables either vanishes
everywhere or on at most 3 lines through the
origin. From this observation it is easy to see that:
\begin{lemma}
Every real cubic polynomial in two variables
can be brought to precisely one of the 5
normal forms in table~\ref{tbl:twovarcubics} by linear change
of variables. The second column provides
the isotropy group of the normal form. The third
describes the condition on a cubic form under
which it has the given normal form.
where $\Sigma_3$, the group of permutations of 3 letters,
is represented by
\[
\Sigma_3 =
\left \{
\begin{pmatrix}
1 & 0 \\
0 & 1
\end{pmatrix},
\begin{pmatrix}
0 & -1 \\
1 & -1
\end{pmatrix},
\begin{pmatrix}
0 & -1 \\
-1 & 0
\end{pmatrix},
\begin{pmatrix}
-1 & 1 \\
0 & 1
\end{pmatrix},
\begin{pmatrix}
1 & 0 \\
1 & -1
\end{pmatrix},
\begin{pmatrix}
-1 & 1 \\
-1 & 0
\end{pmatrix}
\right \}.
\]
\end{lemma}
\begin{proof}
Projectivizing, the points $[x:y]$ in $\mathbb{CP}^1$
on which the cubic form vanishes must be a triple of 
complex points, unless the cubic form vanishes. Because the cubic form is real,
the triple of points must be invariant under complex conjugation.
There could be three distinct real points, or 
a real double point and a distinct real point,
or a triple real point, or a real point and a pair
of complex conjugate points. For instance, if there
is a real point and a pair of conjugate complex 
points, we can arrange that the real point move
to $[x:y]=[1:-1]$, because the group of 
linear transformations acts triply transitively
on $\mathbb{RP}^1$. Lets look at the two complex
points now. The subset $\mathbb{RP}^1 \subset \mathbb{CP}^1$
splits $\mathbb{CP}^1$ into two hemispheres. Picking one of
our complex points, we can
check with a little calculation that the real
linear fractional transformations which fix
it and fix $[1:-1]$ are a finite group. Therefore
the real linear fractional transformations fixing
$[1:-1]$ have a 2-dimensional orbit on 
the nonreal points of $\mathbb{CP}^1$, so we
can put our complex point wherever we like,
as long as we keep it in the same hemisphere.
So we put it at $\left[1/2+i \sqrt{3}/2:1\right]$.
\end{proof}

Of course the dimension of the manifold of cubic
forms having a given normal form is the codimension
of the isotropy group. Thus the generic cubic form
has normal form among the last two on our list.
Indeed these are precisely the nonsingular cubic
forms, which we see immediately from the
discriminant (see \cite{Hilb}, pg. 56), which is
\[\delta(ax^3 + 3bx^2y + 3cxy^2 + dy^3)
= a^2d^2 - 3 b^2c^2 + 4b^3d + 4ac^3 - 6abcd\]
and is positive on $\frac{1}{6} x^3+\frac{1}{6} y^3$
and negative on $\frac{1}{2} x^2 y-\frac{1}{2} x y^2$.

\begin{proposition} The stable cubic forms in two variables
are those on which
the discriminant does not vanish. Everything else
is unstable. The moduli space $\mathfrak M_{3,2}$ is a
pair of points.
\end{proposition}

\section{Cubic Forms in 3 Variables}

\begin{table} 
  \begin{tabular}{|l|l|l|} \hline
    \textbf{Normal form}       
  & \textbf{Isotropy group} 
  & \textbf{Type} \\ \hline
        \fixUp{0} &
        \fixUp{GL(3,\mathbb R)} &
        \\ \hline
        \fixUp{\frac{x^3}6} & 
        \fixUp{\begin{pmatrix}
                1 & 0 & 0 \\
                \cdot & \cdot & \cdot \\
                \cdot & \cdot & \cdot 
        \end{pmatrix}}
        &
        perfect cube
        \\ \hline
        \fixUp{\frac{x^2y}2} &
        \fixUp{\begin{pmatrix}
                a & 0 & 0 \\
                0 & \frac{1}{a^2} & 0 \\
                \cdot & \cdot & \cdot 
        \end{pmatrix}}
        & \fixUp{\text{square} \cdot \text{ind linear}}
        \\ \hline
        \fixUp{\frac{x^2y}2 - \frac{xy^2}2} &
        \fixUp{\begin{pmatrix}
                \Sigma_3 & {} & 0 \\
                {} & {}& 0 \\
                \cdot & \cdot & \cdot   
        \end{pmatrix}}
        & 3 distinct lin dep factors
        \\ \hline
        
        \fixUp{\frac{ x(x^2 + y^2) }{2}} &
        \fixUp{
        \begin{pmatrix}
                1 & 0 & 0 \\
                0 & \pm 1 & 0 \\
                \cdot &  \cdot & \cdot \\
        \end{pmatrix}}
        & \fixUp{\text{lin} \cdot \text{dep semidef quad}}
        \\ \hline
        \fixUp{xyz} & \fixUp{P \cdot
        \begin{pmatrix}
        a & 0 & 0 \\
        0 & b & 0 \\
        0 & 0 & (ab)^{-1}
        \end{pmatrix}}
        & 3 ind factors
        \\ \hline
        \fixUp{
        \frac{ z(x^2 + y^2) }{2}} &
        \fixUp{
        \begin{pmatrix}
                a \cdot O(2) & & 0 \\
                 & & 0 \\
                0 & 0 & a^{-2}
        \end{pmatrix}}
        & \fixUp{\text{lin} \cdot \text{ind semidef quad}}
        \\ \hline
        \fixUp{\frac{x(xz - y^2)}{2}} &
        \fixUp{       
        \begin{pmatrix}
                a^2 & 0 & 0 \\
                b & a^{-1} & 0 \\
                b^2a^{-2} & 2ba^{-3} & a^{-4} 
        \end{pmatrix}}
        & \fixUp{\text{null lin} \cdot \text{lorentz quad}}
        \\ \hline
        \fixUp{\frac{z(x^2 + y^2 - z^2)}{2}} &
        \fixUp{O(2)_z}
        & \fixUp{\text{def lin} \cdot \text{lorentz quad}}
        \\ \hline
        \fixUp{
        \frac{x(x^2 + y^2 - z^2)}{2}} &
        \fixUp{
        \begin{pmatrix}
                1 & 0 & 0 \\
                0 & a & b \\
                0 & b & a \\
        \end{pmatrix}
        \begin{pmatrix}
                1 & 0 & 0 \\
                0 & \pm 1 & 0 \\
                0 & 0 & \pm 1 \\
        \end{pmatrix}, a^2 - b^2 = 1}
        & \fixUp{\text{split lin} \cdot \text{lorentz quad}}
        \\ \hline
        \fixUp{\frac{ x(x^2 + y^2 + z^2) }{2}} &
        \fixUp{O(2)_x}
        & \fixUp{\text{lin} \cdot \text{def quad}}
        \\ \hline
        \fixUp{
        \frac{x^3}6 - \frac{y^2z}2} &
        \fixUp{
        \begin{pmatrix}
                1 & 0 & 0 \\
                0 & a & 0 \\
                0 & 0 & \frac{1}{a^2}
        \end{pmatrix}}
        & cuspidal              
        \\ \hline
        \fixUp{\frac{x^3}6 + \frac{x^2z}2 - \frac{y^2z}2} &
        \fixUp{y \mapsto \pm y} 
        & real nodal
        \\ \hline
        \fixUp{\frac{x^3}6 - \frac{x^2z}2 - \frac{y^2z}2} &
        \fixUp{y \mapsto \pm y}
        & imag nodal 
        \\ \hline
        \fixUp{
        \frac{x^3}6 + \frac{y^3}6 + \frac{z^3}6 + \sigma xyz, \
        \sigma \ne \frac{-1}{2} } &
        \fixUp{P}
        & nonsingular \\ \hline
  \end{tabular}
\caption{Cubic forms in 3 variables}\label{tbl:threevarcubics}
\end{table}

\begin{lemma} Every real cubic form in 3 variables can be
brought to precisely one of the 15 normal forms described in
table~\ref{tbl:threevarcubics} by linear change of variables. 
The second column provides the isotropy group, the third
the condition on a cubic form that it have the given normal
form.
Here $\Sigma_3$ means the same group of $2 \times 2$ matrices that
occured in the previous lemma, $P$ is the group of permutation
matrices, and $O(2)_x$ is the orthogonal group fixing the
$x$ axis. In the case of a split linear times a
Lorentzian quadratic, $a^2 - b^2=1$. The parameter $\sigma$ is an
invariant of any nonsingular cubic form, and the indicated
isotropy group
is for each fixed value of $\sigma$.
\end{lemma}
\begin{proof} For the last row, see \cite{Cool}. The rest
is easy to calculate. For example, the symmetry group
of $xyz$ must act on projectivized points $[x:y:z] \in \mathbb{RP}^2$
permuting the lines $(x=0), (y=0), (z=0)$. Therefore,
after action of the permutation group $P$, we have
a linear transformation for which all three coordinate
axes are eigenspaces, so given by a diagonal matrix.
Plugging in a diagonal matrix, one immediately finds
the symmetry group as stated for $xyz$.
\end{proof}

\begin{proposition} On the list in the previous lemma, the
 last entry (nonsingular) is stable, and the two next to the 
last entry (nodal) are
semistable. The others are unstable.
Stability corresponds to nonsingularity, while
semistability corresponds to the presence of a single
real double point in the associated cubic curve.
The moduli space $\mathfrak M_{3,3}$ is topologically a
figure eight. The forms corresponding to smooth cubic
curves with one circuit lie on one loop of the figure
eight, while those with two circuits lie on the other.
The point in the middle corresponds to the cubic curves
with double point.
\end{proposition}
\begin{proof} See \cite{MFK}.
\end{proof}

\section{Consequences for Lagrangian Surfaces and 3 folds}

\begin{theorem} A Lagrangian surface which has nonsingular
cubic form invariant at a point, has the same cubic form
invariant at
all nearby points. If its cubic form invariant
is never singular, it has a canonical choice of 
framing, up to switching the legs
when the cubic form is
$$\frac{x^3}6 + \frac{y^3}6$$
or up to the previously indicated $\Sigma_3$ action
when the cubic form is
$$\frac{x^2y}2 - \frac{xy^2}2.$$
In the first [second] case, if the fundamental group has no
subgroup of index 2 [index 2, 3 or 6], then the surface bears two [6] canonical
trivializations.
\end{theorem}
\begin{proof}
The stability of the nonsingular cubic forms arises
because they are defined by the condition
that the discriminant is nonzero, an open condition.
The canonical framing is the choice of linear symplectic
coordinates in which the cubic form acheives normal form.
\end{proof}

The lowest order invariant is the choice of cubic form,
i.e. third order terms of the generating function, which
we have seen belongs to one of five classes, after affine
symplectic transformation.

\begin{theorem} Suppose that $L$ is a Lagrangian surface
in an affine symplectic space, with constant cubic form invariant.
Then if that cubic form is nonsingular, we have the trivializations
indicated in the last result. If the cubic form is
$0$ our Lagrangian surface consists of open subsets of
Lagrangian planes. If the cubic form is $\frac{x^3}6$ then the surface
bears a canonical choice of nonvanishing one form (given
in each tangent space by the function $x$). If the cubic
form
is $\frac{1}{2}x^2y$ then there is a global choice of
a pair of line fields ($x=0, y=0$).
\end{theorem}
\begin{proof}
If the cubic form vanishes, then the locally defined
generating function $S$ satisfies
(in any local affine symplectic coordinates $q,p$)
\[
\pd{^3 S}{q^i q^j q^k}=0,
\]
so that $S$ is a quadratic function, and therefore
the graph of the gradient of $S$ (i.e. the
Lagrangian submanifold) is linear.
If the cubic form is $x^3$, then the symmetry
group preserving $x^3$ in each tangent space,
and also preserves $x$, a linear function
on the tangent space and therefore a 1-form.
\end{proof}

These theorems have the obvious topological consequences: a sphere
is simply connected, so a smooth choice of finitely many
bases of each tangent space is as good as one,
i.e. a trivialization. But the sphere has nontrivial
tangent bundle. Therefore no Lagrangian sphere can have constant
nonsingular cubic form. Moreover, a sphere cannot have cubic
form $\frac{x^3}{6}$, because that would provide a global
nonvanishing 1-form (given in each tangent space by $x$),
so a global nowhere vanishing 1-form, which the sphere does not have.

\begin{theorem} A Lagrangian submanifold in an affine
symplectic space with nonsingular (resp. stable) 
cubic form invariant at a point has a nonsingular (resp. stable) 
cubic form invariant at all nearby points.
\end{theorem}
\begin{proof}
The condition for nonsingularity (or stability)
is an open condition on the 3-jet.
\end{proof}

\begin{theorem} A Lagrangian 3 manifold with everywhere
nonsingular cubic form invariant has a canonical framing
up to permutation. A Lagrangian 3 manifold $L$
with everywhere not unstable cubic form has a canonical homomorphism
$$\pi_1(L) \mapsto \mathbb Z * \mathbb Z.$$
A Lagrangian 3 manifold with constant
cubic form invariant equal to
$$\frac{x^3}6 + \frac{x^2z}2 - \frac{y^2z}2$$
has a canonical choice of framing up to switching the sign of
the second vector in the framing. More generally,
a Lagrangian 3 manifold with constant cubic form invariant
has a connection with holonomy in the
indicated isotropy group.
\end{theorem}
\begin{proof}
This homomorphism $\pi_1(L) \mapsto \mathbb{Z} * \mathbb{Z}$
is merely the monodromy of the canonical framing, since
the moduli space $\mathfrak{M}_{3,3}$ is a figure 8, so
has fundamental group $\pi_1\left(\mathfrak{M}_{3,3}\right)=\mathbb{Z} * \mathbb{Z}$.
\end{proof}

\section{Integrability of the Affine Symplectic Group and Heisenberg's
Canonical Quantization}

\begin{theorem} (Grunewald-van Hove) The polynomials
of degree less than or equal to 2
on an affine symplectic space form a maximal 
``quantizable'' set in the algebra of polynomials.
\end{theorem}
For a proof, and a definition of quantizable (in the sense of
Heisenberg) see \cite{GS}. 

\begin{theorem} The Hamiltonian flow of a function on an
affine symplectic space preserves the affine structure
precisely if it is a polynomial of degree at most 2.
\end{theorem}
\begin{proof} The Hamiltonian vector field in affine
Darboux coordinates,
$$ \dot q = \pd{H}{p} \quad
  \dot p = - \pd{H}{q} $$
is an infinitesimal affine transformation precisely if it
is of the form $constant \ +  \ linear$.
\end{proof}

Suppose that our affine symplectic space is a vector
space $V$ and that the symplectic structure is $\Omega$.
We define an antisymmetric linear map $J : V^* \to V$
by the equation
$$\Omega(J \xi,w) = \xi(w).$$
We can describe the ``Heisenberg Lie algebra'' as the
central $\mathbb R$ extension of the affine symplectic
Lie algebra:
$$
\begin{pmatrix}
JA & J \xi & 0 \\
0 & 0 & 0 \\
0 & 0 & c
\end{pmatrix}
\mapsto
\begin{pmatrix}
JA & J \xi \\
0 & 0 
\end{pmatrix}$$
where $A : V \to V^*$ is symmetric, and $\xi \in V^*$.
Now we identify the matrix on the left hand side with
the Hamiltonian function
$$h(v) = \frac{1}{2} \left < Av , v \right > 
        + \left < \xi , v \right > + c$$
and this is an isomorphism between the Lie bracket
and the Poisson bracket, i.e. if we choose
$$f(v) = \frac{1}{2} \left < Av , v \right > 
        + \left < \xi , v \right > + c$$
$$g(v) = \frac{1}{2} \left < Bv , v \right > 
        + \left < \eta , v \right > + d$$
then
$$\{ f,g \} = \frac{1}{2}  \left < (AJB-BJa)v , v \right > 
        + \left < AJ \xi - BJ \eta , v \right > 
        + \left < \xi , J \eta \right >$$

\begin{lemma}\label{lemma:OneParam} For the generic choice of symmetric
linear map $A : V \to V^*$ on a symplectic vector
space $V$ there exist linear
Darboux coordinates $q_1, p_1, \dots , q_n , p_n$
for $V$ so that $A$ is diagonalized, i.e.
$$A = 
\begin{pmatrix}
a_1     & {}    & {}    & {}    & \dots \\
{}      & b_1   & {}    & {}    & \dots \\
{}      & {}    & a_2   & {}    & \dots \\
{}      & {}    & {}    & b_2   & \dots \\
\vdots  & \vdots & \vdots & \vdots & \ddots 
\end{pmatrix}.
$$
We can also write this
$$ \left < A(q,p), (q,p) \right > 
= \sum \left ( a_j q_j^2 + b_j p_j^2 \right ).$$
\end{lemma}
For a proof of this lemma, see \cite{Will}.

\begin{corollary} The symmetric linear maps $A_j$
$$ \left < A_j(q,p), (q,p) \right > 
= a_j q_j^2 + b_j p_j^2 $$
provide linearly independent Hamiltonian functions
$$h_j(q,p) = \frac{1}{2} \left ( a_j q_j^2 + b_j p_j^2 \right )$$
(harmonic oscillators) which Poisson commute 
with one another, and with the Hamiltonian
function
$$h(q,p) = \frac{1}{2} \left < A(q,p), (q,p) \right > 
= \frac{1}{2} \sum \left ( a_j q_j^2 + b_j p_j^2 \right ).$$
\end{corollary}

\begin{corollary} Suppose that $A$ is also invertible.
(This holds generically.) Define $A_j$ as above, and 
$$\xi_j = A_j A^{-1} \xi.$$
Then the Lie algebra spanned by
$$
\begin{pmatrix}
0       & 0     & 0     \\
0       & 0     & 0     \\
0       & 0     & 1
\end{pmatrix}
,
\begin{pmatrix}
JA_1    & J \xi_1       & 0     \\
0       & 0             & 0     \\
0       & 0             & 0
\end{pmatrix}
, \dots ,
\begin{pmatrix}
JA_n    & J \xi_n       & 0     \\
0       & 0             & 0     \\
0       & 0             & 0
\end{pmatrix}
$$
is commutative, and has dimension $n+1$, and contains
$$\begin{pmatrix}
JA      & J \xi         & 0     \\
0       & 0             & 0     \\
0       & 0             & 0
\end{pmatrix}.$$
Moreover it projects to a commutative n dimensional
subalgebra of the affine symplectic Lie algebra
under
$$\begin{pmatrix}
JB      & J \eta        & 0     \\
0       & 0             & 0     \\
0       & 0             & 0
\end{pmatrix}
\mapsto
\begin{pmatrix}
JB      & J \eta\\
0       & 0     
\end{pmatrix}.$$
\end{corollary}
\begin{proof} The proof is linear algebra, using
the obvious result that $JA_jJA_k = 0$ for $j \ne k$.
\end{proof}

\begin{theorem} Every one parameter subgroup of the
affine symplectic group is contained in an n parameter
abelian subgroup.
\end{theorem}
\begin{proof} By lemma~\ref{lemma:OneParam}, this is true
for generic one parameter subgroups (in the Zariski topology),
hence it follows for any one parameter subgroup.
\end{proof}

Thus we can reinterpret our results very heuristically: 
there is an integrable
system at work, and we are following its action on its
Lagrangian submanifolds. Since these are supposed to ``carry
the physical data'', we want to study the local and
global invariants of Lagrangian submanifolds to see what
sort of invariant data they are capable of carrying.

\section{Extreme Points}

\begin{definition} An (weak) \emph{extreme point} of a submanifold
of an affine space is a point at which the restriction of
an affine function
attains a (degenerate) nondegenerate critical point with
positive (semidefinite) definite Hessian.
\end{definition}

\begin{proposition} For a Lagrangian surface with cubic form invariant
at a given point among 
$0, \frac{1}{6}x^3, \frac{1}{2}x^2y$ this point is a weak extreme point.
If it has cubic form $\frac{1}{6}x^3 + \frac{1}{6}y^3$
then this point is an extreme point. Otherwise it is neither
extreme nor weak extreme.
\end{proposition}
\begin{proof}
For a Lagrangian surface with cubic form invariant $0$,
the generating function vanishes up to fourth order,
say $S(x,y,z) = ax^4+ \dots$. Then the Lagrangian
submanifold is the graph of
\[
p_1 = \pd{S}{x}, p_2 = \pd{S}{y}, p_3 = \pd{S}{z},
\]
all cubic expressions in $x,y,z$, up to higher order
terms. The linear function $p_1$ attains a critical point
when restricted to the Lagrangian submanifold.
\end{proof}
\begin{corollary} A Lagrangian surface with nonsingular
cubic form at all points must be everywhere extreme,
or nowhere
\end{corollary}
\begin{corollary} A compact Lagrangian surface (in an
affine symplectic space) with nonsingular cubic form
at all points must have cubic form
\[ \frac{1}{6}x^3 + \frac{1}{6}y^3 \]
(in some linear coordinates, in each tangent space).
\end{corollary}
\begin{proof} (of the proposition) Consider the problem
in Darboux coordinates, with $p=S'(q)$ as our Lagrangian
surface, and suppose $S(0)=0, S'(0)=0, S''(0)=0$.
Our affine function $f$ can be assumed to take value $0$ 
at $0$. It will vanish on the tangent space, so
it looks like
$$f(q,p) = \sum u_i p_i$$
for $(u_1, u_2) \ne 0$ For $f$ to have a nondegenerate
maximum on $L$ (i.e. on $p=\partial S/\partial q$), 
$$0 = \pd{f}{q_j} = S''(0)u$$
which is vacuously satisfied, and
$$\left ( \frac{\partial^2 S}{\partial q_i \partial q_j} \right )
< 0$$
which says
$$\left ( u_k \frac{\partial^3 S}{\partial q_i \partial q_j 
\partial q_k } \right )
< 0.$$ 
The rest is a simple calculation. For instance if $u_1, u_2 < 0$
then we get a nondegenerate maximum for $f$ in the expected
case $S=\frac{1}{6}x^3 + \frac{1}{6}y^3$.
\end{proof}

\begin{corollary} A Lagrangian submanifold of any dimension
in an affine symplectic space, which has nonsingular 
cubic form at all points, 
must have everywhere weak extreme points.
\end{corollary}
\begin{proof} We can slice by symplectic 3 planes to get 
to the case of a Lagrangian 3 manifold. The resulting
3 manifold clearly must be nonsingular, since a singularity
in its cubic form would provide one in the original cubic
form. Now we use our normal form for nonsingular cubic
forms in 3 variables to get
$$S = \frac{1}{6} x^3 + \frac{1}{6}y^3 + \frac{1}{6}z^3 + \sigma xyz + \dots$$
to 3rd order. Now set $z$ and its conjugate momentum variable to $0$.
\end{proof}

\begin{example} The Clifford torus has
extreme points at all points, since it sits in a sphere.
\end{example}

\section{Differential Equations for Prescribing Cubic Form Invariants} 
\label{sec:diffeq}

Let me fix affine Darboux coordinates on our affine
symplectic space. Then the affine symplectic
group can be identified with the bundle of symplectic
framings:
$$
g=
\begin{pmatrix}
u_1     & \dots         & u_n   & v^1   & \dots         & v^n   & x \\
0       & \dots         & 0     & 0     & \dots         & 0     & 1     
\end{pmatrix} 
$$
where
$$
u_j = \begin{pmatrix}
u_j^1 \\
\vdots \\
u_j^{2n}
\end{pmatrix},
v^j = \begin{pmatrix}
v^j_1 \\
\vdots \\
v^j_{2n}
\end{pmatrix},
x = \begin{pmatrix}
q^1 \\
\vdots \\
q^{2n} \\
p_1 \\
\vdots \\
p_n 
\end{pmatrix}
$$
are vectors in $\mathbb R^{2n}$ and $u$'s and $v$'s form
a symplectic framing: 
\begin{align*}
  \Omega(u_i,u_j) & = \Omega(v^i,v^j) = 0 \\
  \Omega(u_i, v^j) & = \delta_i^j
\end{align*}
where $\Omega = dq^i \wedge dp_i$ is our symplectic form.
We write
\begin{equation*}
g = \begin{pmatrix}
u & v & x \\
0 & 0 & 1
\end{pmatrix}
\end{equation*}
for short. The equation
\begin{equation*}dg =
\begin{pmatrix}
du & dv & dx \\
0 & 0 & 1
\end{pmatrix}
\end{equation*}
describes 1 forms on the affine symplectic group,
while the left invariant one forms are given by, say
$$ g^{-1} dg =
\begin{pmatrix}
\alpha  & \beta         & \omega \\
\gamma  & - {}^t \alpha & \eta \\
0 & 0 & 0
\end{pmatrix}
$$
with $\beta, \gamma$ symmetric
(and I write ${}^tM$ for the transpose of a matrix $M$).
By setting $dg=g \cdot g^{-1} \cdot dg$ we find the
first structure equations:
\begin{align*}
du &= u \cdot \alpha + v \cdot \gamma \\
dv &= u \cdot \beta - v \cdot {}^t \alpha \\
dx &= u \cdot \omega + v \cdot \eta.
\end{align*}
The second structure equations, which we get from
the equation 
$$d(g^{-1}dg) = - (g^{-1}dg) \wedge (g^{-1}dg)$$
are
\begin{align*}
d \alpha &= - \alpha \wedge \alpha - \beta \wedge \gamma \\
d \beta &= - \alpha \wedge \beta + \beta \wedge {}^t \alpha \\
d \gamma &= - \gamma \wedge \alpha + {}^t \alpha \wedge \gamma \\
d \omega &= - \alpha \wedge \omega - \beta \wedge \eta \\
d \eta &= - \gamma \wedge \omega + {}^t \alpha \wedge \eta.
\end{align*}

\begin{lemma} Under the map
$$\begin{pmatrix}
u & v & x \\
0 & 0 & 1
\end{pmatrix}
\mapsto
x
$$
from the affine symplectic group to the affine space,
the symplectic form from the affine space lifts to
$\omega^i \wedge \eta_i$. (N.B. the lifted form will not
be symplectic.)
\end{lemma}
\begin{proof}
The form $\omega^i \wedge \eta_i$ is left invariant, and
the lift of the symplectic form will also be left invariant,
since the group preserves the symplectic form. It suffices
to check that these forms match at the identity element.
\end{proof}

To any Lagrangian submanifold $L$ in the affine symplectic
space, we assign the bundle $B=B_L$, which is the submanifold
of the affine symplectic group given as those elements
$$\begin{pmatrix}
u & v & x \\
0 & 0 & 1
\end{pmatrix}$$
such that $x \in L$ and $u$ a framing for $T_xL$. This is
easily seen to be a principal $P$ bundle, where $P$
is a subgroup of the linear symplectic group:
$$P = \left \{
\begin{pmatrix}
A & 0 \\
0 & {}^t A^{-1}
\end{pmatrix}
\begin{pmatrix}
I & B \\
0 & I
\end{pmatrix}
\right \}.$$
This is the stabilizer of the Lagrangian plane $p=0$,
which we have encountered before. We have the inclusion
of $B \subset ASp$ where $ASp$ is the affine symplectic group.
Let me pull back the forms
$$\alpha, \beta, \gamma, \omega, \eta$$
to $B$.
\begin{lemma} $\eta_1 , \dots , \eta_n$ vanish on $B$
\end{lemma}
\begin{proof} By left invariance, it suffices to show that
these forms vanish on some left translate of $B$, so I can
assume that $B$ contains the identity element of $ASp$,
and show that each $\eta_i$ vanishes on tangent vectors
to $B$ at this point. At the identity element, $\eta_i = dp_i$
and
\begin{equation*}
\begin{pmatrix}
u & v & x \\
0 & 0 & 1 
\end{pmatrix}
= I
\end{equation*}
so $u_i = \frac{\partial}{\partial q^i}.$
\end{proof}

\begin{lemma} $\omega^1 \wedge \dots \wedge \omega^n$ never vanishes
on $B$.
\end{lemma}
The proof is the same as that of the previous lemma.

On our bundle $B$ the second structure equations simplify to
\begin{align*}
d \alpha &= - \alpha \wedge \alpha - \beta \wedge \gamma \\
d \beta &= - \alpha \wedge \beta + \beta \wedge {}^t \alpha \\
d \gamma &= - \gamma \wedge \alpha + {}^t \alpha \wedge \gamma \\
d \omega &= - \alpha \wedge \omega \\
d \eta &= 0 = - \gamma \wedge \omega
\end{align*}
with all the $\omega^j$ independent. The last equation, by Cartan's
lemma, says that there exist unique functions $S_{ijk}$, symmetric
in $j,k$ so that $\gamma = S_{ijk} \omega^k$. Since $\gamma_{ij}=
\gamma_{ji}$ the $S$ functions are symmetric in all subscripts.
Plugging $\gamma_{ij} = S \omega$ into the second structure equations
gives
$$ \left ( dS_{ijk} - S_{ijl} \alpha^l_k - S_{ilk} \alpha^l_j
-S_{ljk} \alpha^l_i \right ) \wedge \omega^k = 0.$$
Again by Cartan's lemma, we then have well defined $DS_{ijkl}$,
the covariant derivatives of the $S$, determined
by
$$DS_{ijkl} \omega^l = dS_{ijk} - S_{ijl} \alpha^l_k - S_{ilk} \alpha^l_j
-S_{ljk} \alpha^l_i.$$
The $DS$ functions are also symmetric in all subscripts.

\section{Interpreting the Equations}

We will need to split our matrix
$$\begin{pmatrix}
u & v & x \\
0 & 0 & 1
\end{pmatrix}$$
into more parts:
\[
\begin{pmatrix}
u & v & x \\
0 & 0 & 1
\end{pmatrix}
=
\begin{pmatrix}
a & b & q \\
c & d & p \\
0 & 0 & 1
\end{pmatrix}
\]
So
\begin{align*}
u_j &= a^k_j \pd{}{q^k} + c^k_j \pd{}{p_k} \\
v^j &= b^k_j \pd{}{q^k} + d^k_j \pd{}{p_k}.
\end{align*}
Suppose that $L$ is a Lagrangian submanifold of $\mathbb R^{2n}$
with the form
$$p = S'(q)$$
and that $S(0)=0, \ S'(0)=0, \ S''(0)=0$.
\begin{lemma} The vector fields
\begin{align*}
u_j &= \pd{}{q^j} + \pd{^2 S}{q^j \partial q^k} \pd{}{p_k} \\
v^j &= \pd{}{p_j}
\end{align*}
form a symplectic framing, with the $u$ vector fields tangent
to the Lagrangian submanifold $p=S'(q)$.
\end{lemma}
This is a minor computation. We can map a point $(q,S'(q))$ 
of our Lagrangian submanifold to the affine symplectic group as
\[
g(q) = \begin{pmatrix}
I & 0 & q \\
S'' & I & S' \\
0 & 0 & 1
\end{pmatrix}
\]
so that the first columns give our vector fields $u$ and $v$.
Through this map we pull down all of our forms on the
affine symplectic group to the Lagrangian submanifold.
\[
dg
=
\begin{pmatrix}
0 & 0 & dq \\
dS'' & 0 & dS' \\
0 & 0 & 0
\end{pmatrix},
\] 
so to find $g^{-1} dg$,
\[
g^{-1} = 
\begin{pmatrix}
I & 0 & -q \\
-S'' & I & -S' + S''q \\
0 & 0 & 1
\end{pmatrix}
\]
which gives
\begin{align*}
g^{-1} dg
&=
\begin{pmatrix}
0 & 0 & dq \\
dS'' & 0 & 0 \\
0 & 0 & 0
\end{pmatrix}
\\ &=
\begin{pmatrix}
\alpha & \beta & \omega \\
\gamma & -{}^t \alpha & \eta \\
0 & 0 & 0
\end{pmatrix}.
\end{align*}
Thus $\gamma = dS''$ or
\[ \gamma_{ij} 
= \pd{^{3} S}{q^i \partial q^j \partial q^k} dq^k
= \pd{^{3} S}{q^i \partial q^j \partial q^k} \omega^k
\]
which shows that the invariants $S_{ijk}$ turn out to be the third
derivatives of the generating function. Their covariant
derivatives, immediate from the above expressions, turn out to
be the fourth derivatives.

\section{Setting Our Invariants to Constants}

  If we have a Lagrangian submanifold on which the cubic
form invariant is constant, we can choose (at least
locally) a symplectic framing on which the $S_{ijk}$
are constant, say $C_{ijk}$. This would be a submanifold of the
affine symplectic group on which the forms
\[
  \gamma_{ij} - C_{ijk} \omega^k
\]
all vanish. With the assumption that $C$ is symmetric
in its indices, a submanifold of the affine symplectic
group on which these forms vanish and on which
\[
  \omega^1 \wedge \dots \wedge \omega^n \ne 0
\]
is precisely such a symplectic framing. We take 
these forms to generate a system of partial
differential equations
\[ 
  \gamma_{ij} - C_{ijk} \omega^k = 0 
\]
whose solution manifolds we will study.

\section{The One Dimensional Case}

  All curves in the symplectic plane are Lagrangian.
Our second structure equations are
\begin{align*}
d \alpha &= - \beta \wedge \gamma \\
d \beta &= -2 \alpha \wedge \beta \\
d \gamma &= -2 \gamma \wedge \alpha \\
d \omega &= - \alpha \wedge \omega - \beta \wedge \eta \\
d \eta &= - \gamma \wedge \omega + \alpha \wedge \eta.
\end{align*}
The tableau for Lagrangian curves is
$$d \eta = - \gamma \wedge \omega \quad \mod \eta.$$
The solutions of course depend on 1 arbitrary function
of 1 variable, since a Lagrangian curve can be seen
as locally the graph of a function.
 
  Setting $\gamma = C \omega$ for some constant $C$, and
adding this to our ideal, we find the tableau:
$$d ( \gamma - C \omega ) 
=
3 C \alpha \wedge \omega$$
So on solutions, $\alpha = A \omega$ for some function $A$.
Again we have 1 function of 1 variable for solutions,
unless $C=0$ in which can we obviously get only lines
as solutions. Thus along any curve in the plane it is 
possible to choice a symplectic framing so that the
function $f$ determined by $\gamma = f \omega$ is
constant.

  Now suppose that we try to set $A$ to a constant.
$$d ( \alpha - A \omega)
=
- C \beta \wedge \omega$$
So again 1 function of 1 variable, and we can arrange
our framing to make $A$ constant along with $f$.
But in this framing $\beta = B \omega$ for some
function $B$. We set this to a constant, and find
$$d ( \beta - B \omega ) = 0.$$
So our system is Frobenius, and there is a finite
parameter family of solutions, given by
\[
g^{-1} dg =
\begin{pmatrix}
A & B & 1 \\
C & -A & 0 \\
0 & 0 & 0
\end{pmatrix}
\omega
\]
which determines our curve:
\[
g(t) = \exp \left ( t M  \right )
\]
(up to reparameterization)
where
\[
M =
\begin{pmatrix} 
A & B & 1 \\
C & -A & 0 \\
0 & 0 & 0
\end{pmatrix}.
\]
By scaling, we can assume that $\det \, M =-1, 0$ or $1$.
If $\det \, M = 1$ then the curve through a point $x$ is given
by
\[
g(t) = 
\begin{pmatrix}
\cos{t} I + \sin{tM} & \left ( \sin{t} I + (\cos{t}-1)M^{-1} \right ) x \\
0 & I
\end{pmatrix}.
\]
The Lagrangian curve is an ellipse. If $det \, M =-1$ then
the curve similarly is a hyperbola, and if $det \, M=0$ the
curve is a parabola. These are the homogeneous Lagrangian
curves in the affine symplectic plane.

\section{Lagrangian Surfaces and 3-Manifolds with Fixed Cubic Form}

\begin{proposition} The real analytic Lagrangian surfaces
with constant cubic form invariant arise in the generality
specified in table~\ref{tbl:surfsols}.
\end{proposition}

\begin{table} 
\[
  \begin{array}{|l|l|l|} \hline
    \mbox{\textbf{Normal form}}        & 
    \mbox{\textbf{Generality}} \\ \hline
    0                         & \mbox{all Lagrangian planes}\\ \hline
    \text{perfect cube}       & \mbox{2 functions of 1 variable} \\ \hline
    \text{linear}_1 \cdot \left ( \text{linear}_2 \right )^2
         & \mbox{3 functions of 1 variable} \\ \hline
    \text{linear} \cdot \left ( \text{irred quadratic} \right )           
         & \mbox{1 function of 2 variables (generic)} \\
\hline
    \text{linear}_1 \cdot \text{linear}_2 \cdot \text{linear}_3       
         & \mbox{1 function of 2 variables (generic)}
  \\ \hline
  \end{array}
\]
\caption{Lagrangian surfaces with fixed cubic form}\label{tbl:surfsols}
\end{table}

\begin{proposition} The real analytic Lagrangian 3 manifolds
with constant cubic form invariant arise in the generality
specified in table~\ref{tbl:threeDsols}.
\end{proposition}
\begin{proof}
The precise meaning of the expression ``$k$ functions
of $n$ variables'' comes from Cartan's work on
systems of partial differential equations,
and is explained in detail in \cite{BCGGG}.
The differential equations for constant
nonsingular cubic invariant are in involution without
any prolongation, except in the case of the Fermat cubic
\[ x^3 + y^3 + z^3 \]
for which a single prolongation reveals involution,
with the same Cartan characters as any of the other nonsingular
cases. Computational technique for handling Cartan's test
is explained in detail in \cite{BCGGG}.
\end{proof}

\newcommand{\mcdot}{\( \cdot \ \)}

\begin{table} 
  \begin{tabular}{|l|l|l|} \hline
    \textbf{Normal form}        & 
    \textbf{Generality} \\ \hline
    $0$                            & all Lagrangian affine planes\\ \hline
    perfect cube                   & 3 functions of 1 variable \\ \hline
    square \mcdot
             indep linear          & 5 functions of 1 variable \\ \hline
    3 distinct lin dep factors     & 1 function of 2 variables \\ \hline
    lin \mcdot dep semidef quadratic & 1 function of 2 variables \\ \hline
    3 indep lin factors            & 1 constant \\ \hline
    linear \mcdot ind semidef quad     & 1 constant \\ \hline
    null linear \mcdot lorentz quad   & 1 function of 2 variables \\ \hline
    def linear \mcdot lorentz quad & 8 functions of 1 variable \\ \hline
    split linear \mcdot lorentz quad & 8 functions of 1 variable \\ \hline
    linear \mcdot def quad              & 2 functions of 2 variables \\ \hline
    cuspidal                       & 2 functions of 2 variables \\ \hline
    imaginary nodal                & 3 functions of 2 variables \\ \hline
    real nodal                     & 3 functions of 2 variables \\ \hline
    nonsingular                    & 3 functions of 2 variables \\ \hline  
  \end{tabular}
\caption{Lagrangian 3 manifolds with fixed cubic form}\label{tbl:threeDsols}
\end{table}

\section{Ruled Lagrangian Submanifolds}

In the notation established in section~\ref{sec:diffeq}, the tangent
space to a Lagrangian submanifold in an adapted frame is
\[ u_1 , \dots , u_n. \]
Moving along the adapted frame bundle, we find
\begin{align*}
 du_i &= u_j \alpha^j_i + v^j \gamma_{ji} \\
      &= u_j \alpha^j_i + v^j S_{ijk} \omega^k
\end{align*}
and
\begin{align*}
  d \left ( u_1 \wedge \dots \wedge u_n \right )
  =& \sum_i \left ( -1 \right )^{i+1} \alpha^i_i u_1 \wedge \dots \wedge u_n
\\+& \sum_i \left ( -1 \right )^{i+1} u_1 \wedge \dots \wedge u_{i-1}
  \wedge S_{ijk} \omega^k v^j \wedge u_{i+1} \wedge \dots \wedge u_n.
\end{align*}
We see that the tangent space is carried along in parallel by a vector
tangent to our Lagrangian submanifold precisely if that vector is
in the plane field determined by $\gamma=0$, i.e. by 
$S_{ijk} \omega^k = 0$. These are precisely the tangent directions
in which our cubic form on the given tangent space is invariant under translation, i.e. the ``unused variables'' of the cubic form.
For instance in $\mathbb R^n$, the form $C(x_1 , \dots , x_k)$
has among its ``unused variables'' any linear combination of 
\[ 
  \pd{}{x^{k+1}}, \dots, \pd{}{x^n}.
\]
This plane field \( \gamma = 0 \) is holonomic,
since our structure equations give 
\[
d \gamma = - \gamma \wedge \alpha + {}^t \alpha \wedge \gamma.
\]
Moreover, it is a field of Cauchy characteristics
(see \cite{BCGGG}), so that we can reduce the differential
system.

  Let \( \AffLag{d}{n} \) be the space of affine sublagrangian
$d$ planes in $\mathbb R^{2n}$. If we let \( \Lag{d}{n} \)
be the linear sublagrangian $d$ planes, then 
\[
\AffLag{d}{n} \to \Lag{d}{n}
\]
is a vector bundle of rank $2n-d$. The dimensions are
\begin{align*}
\dim \Lag{d}{n} &= k \frac{4n-3k + 1}{2} \\
\dim \AffLag{d}{n} &= k \frac{4n-3k + 1}{2} + 2n - k.
\end{align*}
We can include \( \Lag{d}{n} \subset Gr_d ( \mathbb R^{2n} ) \),
and identify
\[
T_P Gr_d ( \mathbb R^{2n} ) \cong \mbox{Lin}(P, \mathbb R^{2n}/P)
\]
invariantly under the $GL(2n, \mathbb R)$ action. The
subspace of tangent vectors to \( \Lag{d}{n} \) is identified
as follows: given a tangent vector \( 
A \in \mbox{Lin}(P, \mathbb R^{2n}/P) \)
we can define the bilinear form 
\[ v,w \in P \mapsto \Omega ( A v, w ) \]
where $\Omega$ is our symplectic form on $\mathbb R^{2n}$.
\[
T_P \Lag{d}{n} \cong \{ A \in \mbox{Lin}(P, \mathbb R^{2n}/P) :
\Omega(A \cdot , \cdot ) \mbox{ antisymmetric} \}
\]
This gives a well defined section $\sigma$ of 
the bundle \( T^* \Lag{d}{n} \otimes \Lambda^2 U^* \)
where $U$ is the universal bundle on $\Lag{d}{n}$
by
\[
\sigma(A,v,w) = \Omega(Av,w).
\]

I will need to refine my notation a little to explain 
how the Pfaffian system on $\AffLag{d}{n}$ is
determined. I will write my left invariant one forms
on the affine symplectic group splitting them into
subscripts running $1, \dots ,d$ and $d+1, \dots, n$,
as, for example
\begin{align*}
\alpha &= 
\begin{pmatrix}
( \alpha )_{11} & ( \alpha )_{12} \\
( \alpha )_{21} & ( \alpha )_{22}
\end{pmatrix} \\
\omega &=
\begin{pmatrix}
( \omega )_1 \\
( \omega )_2
\end{pmatrix}.
\end{align*}
This unfortunate notation allows me to distinguish between
$\alpha^2_2$ and $( \alpha )_{22}$, which hopefully the reader
can distinguish as well.

The tangent spaces to the fibers of the projection
\begin{align*}
\pi  : \ASp{n} & \to \AffLag{d}{n} \\
       \begin{pmatrix}
        u & v & x \\
        0 & 0 & 1
       \end{pmatrix} & \mapsto d \mbox{ plane through x containing }
u_1 , \dots , u_d
\end{align*}
are given by the equations
\[
( \alpha )_{21} = ( \gamma )_{11} = ( \gamma )_{12} = \eta = ( \omega )_2 = 0.
\]
These equations cut out the Lie algebra of the isotropy group of
\( \mathbb R^d \subset \mathbb R^{2n} \) among the left
invariant vector fields.

   Now if we pick a cubic form \( C \in S^3 \mathbb R^n \) with
its $d$ unused variables being $x^1 , \dots , x^d$, then we 
can consider the Pfaffian system \( \gamma = C \cdot \omega, \eta = 0 \)
on \( \ASp{n} \) and try to ``push'' it down to \( \AffLag{d}{n} \)
along the projection. On each fiber, \( \ASp{n}_P \) 
of \( \pi  : \ASp{n} \to \AffLag{d}{n} \)
we have a foliation
given by the holonomic plane field \( \gamma_{ij} = C_{ijk} \omega^k \).
The problem of solving the Pfaffian system we started with on
\( \ASp{n} \) should resolve in to two problems: solving the
reduced system on \( \AffLag{d}{n} \), and then choosing the
correct leaf of this foliation above each $d$ plane.

\section{The Tableau Product}
We will make an algebraic description of the space
of integral elements of our Pfaffian equations.
\begin{definition} We define a product
$$A \otimes C \in Lin(V \otimes V, V) \otimes S^3V^*
\mapsto
AC \in S^2V^* \otimes \Lambda^2V^*$$
given by
\begin{align*}
AC(u_1,u_2,u_3,u_4)
=&
C(u_1, u_2, A(u_3, u_4))
+
C(u_1, A(u_2, u_4), u_3)
+
C(A(u_1, u_4), u_2, u_3) \\
&
-
C(u_1, u_2, A(u_4, u_3))
-
C(u_1, A(u_2, u_3), u_4)
-
C(A(u_1, u_3), u_2, u_4)).
\end{align*}
We define another product
\[
B \otimes C \in S^2V \otimes S^3V^*
\mapsto
BC \in Lin(V \otimes V , V)
\]
by
\[
\left < \lambda , BC(u,v) \right > := C(u,v,B \lambda)
\]
where we take $u,v \in V, \lambda \in V^*$ and think of
\( B \) as a symmetric linear map \( B : V^* \to V \).
\end{definition}

\begin{proposition} For any \( C \in S^3V^* \), \( C \) is
singular at \( w \in V \) iff \( AC=0 \) for some
(and hence any)
$A \in Lin(V \otimes V , V)$ of the form
$$A(u,v)=w \left < L(u), v \right >$$
with $L : V \to V^*$ any nonzero linear map.
Consequently, for $C$ singular, the kernel
of $$A \mapsto AC$$ in $Lin(V \otimes V,V)$ 
has dimension at least $n^2$. This dimension
is $n^3$ precisely for $C=0$.
\end{proposition}

\begin{proposition} If $A = BC$ then $AC = 0$. The
map $B \mapsto BC$ is injective from $S^2V$ for
generic $C \in S^3V^*$, hence
$$A \mapsto AC$$
has kernel in $Lin(V \otimes V , V)$ of dimension
at least $n(n+1)/2$.
\end{proposition}

\begin{proposition} For any $C \in S^3V^*$, the integral
elements of the equations given
above are parameterized by the vector space
$$\{ A \in Lin(V \otimes V , V) : AC = 0 \}.$$
\end{proposition}

Therefore to apply the Cartan-K\"ahler theorem to our
differential equations, we have only to
count the dimension of this vector space,
and compare to the sum $\sum ks'_k$ given
by the reduced Cartan characters $s'_k$.

\begin{example} Consider the cubic form
$$C(x)= \frac{1}{6}\left < \lambda , x \right >^3$$
for some nonzero $\lambda \in V^*$. 
We can choose linear coordinates $x^1, \dots , x^n$
on our vector space in which $\lambda = x^1$. Then
it is easy to compute that $AC=0$ precisely
if
$$A^1_{ij} = 0 \ i,j>1 \quad A^1_{j1} = 3A^1_{1j}.$$
Therefore the kernel of $A \mapsto AC$ has dimension
$n^2(n-1)+n$. 
\end{example}

\end{document}